\documentclass[11pt]{article}
\usepackage{authblk}
\usepackage[T1]{fontenc}
\usepackage[utf8]{inputenc}
\usepackage[english]{babel}
\usepackage{lmodern}
\usepackage{mathtools} 
\usepackage[a4paper,DIVcalc]{typearea} 
\typearea[10mm]{16} 

\providecommand{\abstractt}[1]
{{
  \small
  \textbf{Abstract -} #1
}}

\providecommand{\keywords}[1]
{{
  \small 
  \textbf{Keywords -} #1
}}

\providecommand{\AMSClassification}[1]
{{
  \small
  \textbf{AMS Subject classifications -} #1
}}

 \usepackage{comment}
\usepackage{cases}
\usepackage{setspace}
\usepackage{amsthm,amssymb,amsmath,eqnarray}
\usepackage{dsfont}
\usepackage{latexsym}
\usepackage[pdftex]{graphicx}
\usepackage{mathdots}
\usepackage[hyperindex=true, pdftex=true, linkcolor=black,colorlinks=true,citecolor=blue,urlcolor=blue]{hyperref}
\usepackage{amssymb,amsfonts,textcomp}
\usepackage{multirow}
\usepackage{longtable}
\usepackage{xcolor}
\definecolor{LightGray}{gray}{0.9}

\usepackage{pgf, tikz}
\usepackage{authblk}

\title{Admitted symmetries of Backward Stochastic Differential Equations  }
\author{OUKNINE Anas, LESCOT Paul}

\affil{Laboratoire de Math\'ematiques Rapha\"el Salem, University of Rouen, UMR CNRS 6085, Avenue de l'Universit\'e, 76801 Saint Etienne du Rouvray, France\\
Laboratoire de Math\'ematiques Rapha\"el Salem, University of Rouen, UMR CNRS 6085, Avenue de l'Universit\'e, 76801 Saint Etienne du Rouvray, France
\\
Email: \texttt{anas.ouknine@univ-lemans.fr}\\
Email: \texttt{paul.lescot@univ-rouen.fr}}
\date{\today}

\newtheorem{thm}{Theorem}[section]
\newtheorem{defin}[thm]{Definition}
\newtheorem{prop}[thm]{Proposition}

\newtheorem{remark}[thm]{Remark}

\numberwithin{equation}{section}
\theoremstyle{remark}

\usepackage{cases}
\usepackage{empheq}
\usepackage{hyperref}

 \allowdisplaybreaks[1]

\begin{document}

\maketitle


\noindent \abstractt{In this article, we introduce the concept of admitted Lie group of transformations for both backward stochastic differential equations (BSDEs) and forward backward stochastic differential equations (FBSDEs), following the approach of Meleshko et al. An application to BSDE is presented.}
\newline

\noindent \keywords{Stochastic process; Backward stochastic differential equation; Lie group of transformations; Random time change;  Forward–backward
stochastic differential equation}\newline

\noindent\AMSClassification{34A26, 91B28, 60H10, 60H30, 58D19}

\section{Introduction}
The Lie symmetry method is a elegant and powerful tool to solve differential equations, including both ODEs and PDEs.
By uncovering the symmetries of differential equations, this approach often enables one to derive explicit, and sometimes very intricate, solutions. Such solutions not only provide a direct answer to the differential equation but also enable a deeper analysis of the equation's properties. For example, exploiting these symmetries can lead to a reduction in the order of the equation or a transformation into a simpler form that is easier to manage and analyze. A general survey of symmetry groups for differential equations can be found in \cite{OLV86,ovsiannikov2014group,ibragimov1995crc,blu, hydon2000symmetry,olver1995equivalence,lescot2019symmetries,LESZA,anas2024lie}.
the main idea behind the Lie symmetry method is to analyze the local groups of transformations that map solutions of a differential equation to other solutions of the same equation.\\
In contrast, the application of symmetry techniques to stochastic differential equations (SDEs) has received relatively little attention. The stochastic setting introduces additional challenges due to randomness, and as a result, only a handful of studies have explored symmetry methods in this context \cite{thieullen1997symmetries,huang2023second,srihirun2007definition,meleshko2010new,gaeta2017random,Giuseppe-Gaeta_1999,Giuseppe-Gaeta_2000,gaeta2019w, unal2003symmetries,kozlov2010symmetries,kozlov2011lie,kozlov2020symmetries,kozlov2021symmetries,nass2016n,nass2016symmetry,nass2017w}. Specifically, only three articles investigated the symmetries of backward stochastic differential equations (BSDEs)\cite{sym11091153,ZHANG2021105527,ouknine2025symmetry}.\\
In 1973, Bismut \cite{bismut1973conjugate} was solving a stochastic control problem, and got a non linear BSDE. Later in 1990 Pardoux-Peng \cite{PARDOUX199055} generalized the linear BSDEs that appears in Bismut's work. Under suitable conditions, including regularity and Lipschitz assumptions on the generator, they proved the existence and uniqueness of solutions. Their work was a major breakthrough, forming the foundation for further advances in stochastic analysis and mathematical finance and beyond.
In summary, while the Lie symmetry method has long been a cornerstone in the anlysis of deterministic differential equations, its extension to the stochastic realm an active and promising area of research, with significant contributions already shaping our understanding of BSDEs and related problems. 

Consider a BSDE,
\begin{align}
    \begin{cases}
       dy_t &= -g(t, y_t, z_t) dt + z_t dB_t,\\\label{BSDE*}
         y_{_T} &= \xi,
    \end{cases}
\end{align}
where $T>0$ is the terminal time, and $B_t$ is a standard Brownian motion on $\mathbb{R}$ defined on a probability space $(\Omega,(\mathcal{F}_t)_{0\leq t\leq T},\mathcal{F},\mathbb{P})$, and $(\mathcal{F}_t)_{0\leq t\leq T}$ is the filtration generated by $B_s$ and $s\leq t$. The generator $g$ is defined on $[0,T]\times \mathbb{R}\times \mathbb{R}$ and is assumed to always be measurable with respect to $\mathcal{B}([0,T])\times \mathcal{B}(\mathbb{R})\times \mathcal{B}(\mathbb{R})$. Additionally, for any $y$ and $z$, the process $t\longrightarrow g(t,y_t,z_t)$ is progressively measurable. 
The terminal condition $y_T$ is a $\mathcal{F}_T$-measurable random variable, typically denoted by $\xi$. The solution $(y_t,z_t)_{t\in [0,T]}$ consists of two stochastic processes $(y_t)_{t\in [0,T]}$ and $(z_t)_{t\in [0,T]}$ which are adapted to the filtration $(\mathcal{F}_t)_{t\in [0,T]}$. \\
Since then, BSDE theory has rapidly evolved, both in theoretical advancements and practical applications, such as finance and stochastic control. However, obtaining exact solutions for BSDEs remains a challenging task, with few studies addressing this issue. The Lie-point symmetries method, due to its property, can be effectively applied to BSDEs under any conditions, provided a solution exists.

The first attempt to compute Lie symmetries of BSDEs was made by 
Zhang and Jia \cite{sym11091153} specifically focusing on fiber preserving transformations 
\begin{equation}
    \bar{y}=\phi(t,y,\varepsilon),\quad \bar{t}=\varphi(t,\varepsilon).
\end{equation}
The second approach focuses on W-symmetries/random symmetries of BSDE, involving transformations of the form
\begin{equation}
    \bar{y}=\psi(t,y,w,\varepsilon),\quad \bar{t}=\varphi(t,\varepsilon), \quad \bar{w}=h(y,t,w,\varepsilon).
\end{equation}
where $\psi$  and $h$ depend on the Brownian motion also. W-symmetries were introduced for the first time by Gaeta \cite{Giuseppe-Gaeta_2000}, who applied them to stochastic differential equations (SDEs).

This transformation using Ito formula maps \eqref{BSDE*} into the equation 
\begin{equation}
    d\bar{y}_{\bar{t}} = -g(\bar{t}, \bar{y}_{\bar{t}}, \bar{z}_{\bar{t}}) d\bar{t} + \bar{z}_{\bar{t}} d\bar{B}_{\bar{t}}
\end{equation}
Using Ito formula \cite{oksendal2013stochastic,revuz2013continuous}, the evolution of a scalar function $I(t,y)$ satisfies the condition 
\begin{equation}
dI=I_ydY+I_tdt+\frac{1}{2}I_{yy}d<Y,Y>.
\end{equation}
An admitted symmetry of a BSDE \eqref{BSDE*} is an infinitesimal transformation that maps any solution of \eqref{BSDE*} to another solution of the same equation. This approach was applied to SDEs and BSDEs \cite{sym11091153,Giuseppe-Gaeta_1999} but only to fiber-preserving symmetries. The second approach \cite{unal2003symmetries,ibragimov2004approximate} deals with symmetry transformations for SDE involving all the dependent variables in the transformation. In this
approach there is also an attempt to involve Brownian motion in the transformation. Unfortunately, there is
no strict proof that the transformed Brownian motion remains a Brownian motion.\\
When analyzing the admitted symmetries of a BSDE \eqref{BSDE*}, it is obvious to consider the transformation group associated with the spatial variable $y$ and the temporal variable $t$. However, the process $z$, which appears in BSDE solutions, raises an important question: should it be explicitly included in the transformation ? The answer is no. Since $z_t$
inherently accompanies $y_t$, its transformation can be derived automatically once the transformation for $y$ is determined. 
This manuscript focuses on deriving the determining equations for admitted Lie groups of transformations, incorporating both the dependent variables and the Brownian motion within the transformations. This approach has already been developed for SDEs \cite{srihirun2007definition,srihirun2006definition,meleshko2010new}, but not for BSDEs.
\section{Random time change}
Let $(\Omega,(\mathcal{F}_t)_{_{0\leq t\leq T}},\mathcal{F},\mathbb{P})$ a filtered probability space. The following constructions are analogous to \cite{srihirun2007definition,srihirun2006definition,oksendal2013stochastic,oksendal1985stochastic}. Let $\eta(t,y,a)$ be a sufficiently many times continuously differentiable function and $(Y_t)_{t\geq 0}$ a continous and $\mathcal{F}_t$-adapted stochastic process. Since $\eta^2(t,y,a)$ is continuous, $\eta^2(t,Y(t,\omega),a)$ is also an $\mathcal{F}_t$-adapted process. Define
\begin{equation}
    \beta(t,\omega,a)=\int_0^t \eta^2(t,Y(t,\omega),a)ds,\quad t\geq 0.\label{rate beta}
\end{equation}
For brevity, we denote $\beta(t)$ instead of $\beta(t,\omega,a)$. The function $\beta(t)$ is called a random time change with time change rate $\eta^2(t,Y(t,\omega),a)$. Note that $\beta(t)$ is an $\mathcal{F}_t$-adapted process. Suppose now that $\eta(t,y,a)\neq 0$ for all $(t,y,a)$. Then, the map $t \longrightarrow \beta(t)$ is strictly increasing. Next define
\begin{equation}
    \alpha(t,\omega,a)=\inf\{s:\beta(s,\omega,a)>t\},\label{rate alpha}
\end{equation}
we write also $\alpha(t)$ instead of $\alpha(t,\omega,a)$. for each $\omega$, the map $t \longrightarrow \alpha(t)$ in nondecreasing and continous. One can show that for almost all $\omega$, and for all $t\geq 0$, 
\begin{equation}
    \beta(\alpha(t))=\alpha(\beta(t))=t.
\end{equation}
since $\beta(t)$ is $\mathcal{F}_t$-adapted process, one has $\{\omega: \alpha(t)\leq s\}=\{\omega: t\leq \beta(s)\}\in \mathcal{F}_s$ for all $t\geq 0$ and $s\geq 0.$ Hence $t\longrightarrow \alpha(t)$ is an $\mathcal{F}_s$- stopping time for each $t$.
The theorem below will be essential for defining the transformation of a Brownian motion.
\begin{thm}\cite{srihirun2007definition}\label{new bar(B) Meleshko}
    Let $\eta(t,y,a)$ and $Y_t$ be as above and $B_t$ a standard Brownian motion. Define 
    $$\Bar{B}_t=\int_0^t\eta(s,Y(s,\omega),a)dB_s,\quad\quad t\geq 0$$
    Then $(\Bar{B}_{\alpha(t)},\mathcal{F}_{\alpha(t)})$ is a standard Brownian motion where \\
$\mathcal{F}_{\alpha(t)}=\{A\in\mathcal{F}:A\cap\{\omega:\alpha(t)\leq s\}\in\mathcal{F}_s\quad \text{for all}\,\, s\geq 0\}$ and $\alpha$ as defined above.
\end{thm}
\noindent For a full proof of the theorem, 
 we refer to \cite{srihirun2007definition}.\\
\section{Lie group of transformations for stochastic process}
Suppose that the set of transformations
\begin{align}
    \bar{y}&=\varphi(t,y,a)\label{lie group meleshko for Y(BSDE)}\\
    \bar{t}&=\Xi(t,y,a)\label{lie group meleshko for time t(BSDE)}
\end{align}
composes a Lie group.
These transformations are related (generated) by the the infinitesimal generator (in the sense of Lie groups)  
\begin{equation}
    h(t,y)\partial_t+\gamma(t,y)\partial_y \label{vector field BSDE meleshko}
\end{equation}
 \begin{equation}
  \text{such that},\quad  h(t,y)=\frac{\partial \Xi(t,y,a)}{\partial a}|_{a=0}=\Xi_a(t,y,0),\quad \gamma(t,y)=\frac{\partial \varphi(t,y,a)}{\partial a}|_{a=0}=\varphi_a(t,y,0)\label{h&xi-Meleshko(BSDE)}
\end{equation} 
and $\Xi$, $\varphi$ satisfy Lie's equations \cite{OLV86}
\begin{equation}
    \frac{\partial \Xi}{\partial a}=h(\Xi,\varphi),\quad  \frac{\partial \varphi}{\partial a}=\gamma(\Xi,\varphi)\label{Lies equation for BSDE}
\end{equation}
with initial conditions for $a=0$,
\begin{equation}
    \Xi(t,y,0)=t,\quad \varphi(t,y,0)=y.\label{initial cond for lie equ BSDE}
\end{equation}
Since \( \Xi_t(t,y, 0) = 1 \), then \( \Xi_t(t,y,a) > 0 \) in a neighborhood of \( a = 0 \), where one can find a function \( \eta(t,x,a) \) such that:
\begin{equation}
\eta^2(t, y, a) = \Xi_t(t, y, a).
\end{equation}
Using the function \( \eta(t,y,a) \), one can define a transformation \( \bar{Y}(\bar{t},\omega) \) of a stochastic process \( Y(t,\omega) \) by:
\begin{equation}
 \bar{Y}(\bar{t},\omega)  = \varphi(\alpha(\bar{t}), Y(\alpha(\bar{t}),\omega), a),\label{bar de Y en fct de bar de t et omega}
\end{equation}
where the functions \( \beta(t) \) and \( \alpha(\bar{t}) \) are as in formulae \eqref{rate beta} and \eqref{rate alpha}. This gives the action of Lie group transformations \eqref{lie group meleshko for Y(BSDE)} and \eqref{lie group meleshko for time t(BSDE)} on the set of stochastic processes. replace $\bar{t}$ by $\beta(t)$, 
\begin{equation}
    \bar{Y}(\beta(t),\omega)=\varphi(t,Y(t,\omega),a).
\end{equation}
For calculations of an admitted Lie group of transformations, it is useful to introduce the function:
\begin{equation}
\tau(t, y) = \eta_a(t, y, 0).
\end{equation}
The functions \( h(t,y) \) and \( \tau(t,x) \) are related by the formulae:
\begin{equation}
\tau(t, y) = \frac{h_t(t, y)}{2}, \quad h(t, y) = 2 \int_0^t \tau(s, y) ds.
\end{equation}

Similar to PDEs, the functions \( \tau(t,y) \) and \( \gamma(t,y) \) define a Lie group of transformations for stochastic processes. in fact, given \( \tau(s,y) \) and \( \gamma(s,y) \), one sets 
$$h(t, y) = 2 \int_0^t \tau(s, y) ds,\quad \text{with} \quad h(0,y)=0.$$
Solving the Lie equations \eqref{Lies equation for BSDE} with initial conditions \eqref{initial cond for lie equ BSDE}, one then finds the functions $\Xi(t,y,a)$ and $\varphi(t,y,a)$.

\section{Admitted symmetries of backward stochastic differential equations}\label{Adm sym of BSDE (paragraph)}

\subsection{Determining equations of BSDE}
This section is dedicated to constructing the determining equations for an admitted Lie group of transformations.\\
Assume that  $(Y(t, \omega),Z(t, \omega))$ are a stochastic process satisfying the following BSDE:
\begin{equation}
    Y(t,\omega) =Y(T,\omega)+\int_{t}^{T}g(s,Y_s,Z_s) ds-\int_{t}^{T}Z_s d B_{s},\quad t\in [0,T],\label{BSDE au sens de Meleshko}
\end{equation}
where the generator $g$ of the BSDE is $\mathcal{F}_t$-adapted. As demonstrated, the process 
 $$\Bar{B}_{\alpha(t)}=\int_0^{\alpha(t)}\eta(s,Y(s,\omega),a)dB(s),\quad\quad t\geq 0$$ is a $\mathcal{F}_{\alpha(t)}$ standard Brownian motion.\\
 Define $$\chi(t,\omega)=\Bar{Y}(\beta(t),\omega),$$ 
 where $\Bar{Y}(t,\omega)=\varphi(\alpha(t),Y(\alpha(t),\omega),a)$ is the transformation of $Y(t,\omega)$ given by \eqref{bar de Y en fct de bar de t et omega}.
 
 For almost all $\omega$,
$$
\chi(t, \omega)=\varphi \left(\alpha(\beta(t)), Y(\alpha(\beta(t)), \omega), a\right)=\varphi(t, Y(t, \omega), a).
$$
By applying Itô's formula, we obtain
\begin{align}
\nonumber \chi(t, \omega)&=\chi(T, \omega)- \int_{t}^{T}
\varphi_{y}(s,Y(s,\omega),a)dY(s,\omega)-\int_{t}^{T} \varphi_{s}(s,Y(s,\omega),a)ds\\
\nonumber&-\frac{1}{2}\int_{t}^{T}
\varphi_{yy}(s,Y(s,\omega),a)d<Y,Y>_s\\
\nonumber&=\varphi(T, Y(T,\omega),a)-\int_{t}^{T}
\varphi_{y}(s,Y(s,\omega),a)(-g(s,Y(s,\omega),Z(s,\omega))ds+Z(s,\omega)dB(s))\\
\nonumber&-\int_{t}^{T} \varphi_{s}(s,Y(s,\omega),a)ds-\frac{1}{2}\int_{t}^{T}
\varphi_{yy}(s,Y(s,\omega),a)Z^2(s,\omega)ds\\
\nonumber&=\varphi(T, Y(T,\omega),a)+\int_{t}^{T}\left(g\varphi_{y}-\varphi_s-\frac{1}{2}\varphi_{yy}Z^2\right)(s,Y(s,\omega),a)ds \label{Ito form BSDE (adm sym)}\\
&-\int_{t}^{T}(\varphi_y Z)(s,Y(s,\omega),a)dB(s)
\end{align}

Because $Y(t, \omega)$ is a solution of \eqref{BSDE au sens de Meleshko} and $\varphi_y(t, y, a)$ is continuous, $\varphi_y(t, Y(t, \omega), a)Z(s,\omega)$ is a bounded and continuous process. Based on the time-change formula for Itô integrals \cite{oksendal2013stochastic}, the bounded and continuous process $v(t, \omega)$ satisfies the formula
\begin{equation}
   \int_{0}^{\alpha(t)} v(s, \omega) \mathrm{d} B(s)=\int_{0}^{t} v(\alpha(s), \omega) \frac{1}{\eta(\alpha(s), Y(\alpha(s), \omega), a)} \mathrm{d} \bar{B}(s)\label{time change formula oksendal}
\end{equation}
Using this formula to the last term of Eq. \eqref{Ito form BSDE (adm sym)}, one finds that
\begin{align}
\nonumber \chi(t, \omega)&=\varphi(T, Y(T,\omega),a)+\int_{\beta(t)}^{\beta(T)}\left(g\varphi_{y}-\varphi_s-\frac{1}{2}\varphi_{yy}Z^2\right)(\alpha(s),Y(\alpha(s),\omega),a)\alpha_{\bar{t}}(s)ds \\
&-\int_{\beta(t)}^{\beta(T)}(\frac{\varphi_y Z}{\eta})(\alpha(s),Y(\alpha(s),\omega),a)dB(s)\label{Ito form BSDE (adm sym)oksendal time change}
\end{align}
Since $\beta(t, \omega, a)=\int_{0}^{t} \eta^{2}(s, Y(s, \omega), a) \mathrm{d} s$ and $\beta(\alpha(\bar{t}))=\bar{t}$ for almost all $\omega$,
$$
\eta^{2}(\alpha(\bar{t}), Y(\alpha(\bar{t}), \omega), a) \alpha_{\bar{t}}(\bar{t})=1
$$
This gives
\begin{equation}
     \alpha_{\bar{t}}(s)=\eta^{-2}(\alpha(s), Y(\alpha(s), \omega), a)\label{alpha_bar(t)(s)}
\end{equation}
Substitution of $\alpha_{\bar{t}}(s)$ into \eqref{Ito form BSDE (adm sym)oksendal time change} gives the equation
\begin{align}
\nonumber \chi(t, \omega)&=\varphi(T, Y(T,\omega),a)+\int_{\beta(t)}^{\beta(T)}\frac{\left(g\varphi_{y}-\varphi_s-\frac{1}{2}\varphi_{yy}Z^2\right)}{\eta^2}(\alpha(s),Y(\alpha(s),\omega),a)\alpha_{\bar{t}}(s)ds \\
&-\int_{\beta(t)}^{\beta(T)}\left(\frac{\varphi_y Z}{\eta}\right)(\alpha(s),Y(\alpha(s),\omega),a)dB(s)
\end{align}

Requiring that transformations \eqref{Lies equation for BSDE} map a solution of the BSDE \eqref{BSDE au sens de Meleshko} into a solution of the same BSDE, one obtains
\begin{equation*}
    \bar{Y}(\bar{t},\omega) =\bar{Y}(\bar{T},\omega)+\int_{\bar{t}}^{\bar{T}}g(s,\bar{Y}(s,\omega),\bar{Z}(s,\omega) ds-\int_{\bar{t}}^{\bar{T}}\bar{Z}(s,\omega) d \bar{B}_{s},\quad t\in [0,T] 
\end{equation*}

Replacing $\bar{t}=\beta(t)$ and $\bar{T}=\beta(T)$ into this equation, one gets
$$
 \bar{Y}(\beta(t),\omega) =\bar{Y}(\beta(T),\omega)+\int_{\beta(t)}^{\beta(T)}g(s,\bar{Y}(s,\omega),\bar{Z}(s,\omega) ds-\int_{\beta(t)}^{\beta(T)}\bar{Z}(s,\omega) d \bar{B}_{s},\quad t\in [0,T] 
$$
Comparing the Riemann integrals and the Itô integrals, we obtain,
\begin{align}
& \frac{\left(g\varphi_{y}-\varphi_s-\frac{1}{2}\varphi_{yy}Z^2\right)}{\eta^{2}}(\alpha(t), Y(\alpha(t), \omega), a)=g(t, \bar{Y}(t, \omega),\bar{Z}(t,\omega))\label{Pre-det-sys-for-BSDE(Mleshko)1}\\
& \bar{Z}(t,\omega)=\left(\frac{\varphi_y Z}{\eta}\right)(\alpha(t), Y(\alpha(t), \omega), a).\label{Pre-det-sys-for-BSDE(Mleshko)2}
\end{align}
We can now give the Lie group that transforms the variable $Z$, which we shall refer to as $\zeta$.
\begin{prop}
     \begin{equation} \zeta(t,Y(t,\omega),Z(t,\omega),a):=\bar{Z}=\frac{\varphi_y Z}{\eta}.\label{bar((Z) meleshko BSDE}
 \end{equation}
\end{prop}
\begin{remark}
In the proposition above, we do not construct a full Lie group structure. Instead, we identify a component of the infinitesimal generator \eqref{vector field BSDE meleshko} corresponding to a one-parameter subgroup of transformations. For this reason, we do not include a term of the form \(\zeta(t,y,z)\,\frac{\partial}{\partial z}\) in \eqref{vector field BSDE meleshko}, as the coefficient \(\zeta\) is computed using the transformation $\varphi$ and the rate $\eta$. 

Moreover, including \(Z\) within the group variables \eqref{vector field BSDE meleshko} would, upon applying Itô's formula, lead to the appearance of a stochastic differential \(dZ\), which is not defined in the current framework. This motivates our choice to restrict the group action to variables for which the stochastic differentials are well-defined.
\end{remark}
Replacing $\eqref{Pre-det-sys-for-BSDE(Mleshko)2}$ in $\eqref{Pre-det-sys-for-BSDE(Mleshko)1}$, we get  
\begin{align}
     \frac{\left(g\varphi_{y}-\varphi_s-\frac{1}{2}\varphi_{yy}Z^2\right)}{\eta^{2}}(\alpha(t), Y(\alpha(t), \omega), a)=g\left(t, \bar{Y}(t, \omega),\left(\frac{\varphi_y Z}{\eta}\right)(\alpha(t), Y(\alpha(t), \omega), a)\right)\label{Pre-det-sys-for-BSDE(Mleshko)2 and 1}
\end{align}
Since $\bar{Y}(\bar{t}, \omega)=\varphi(\alpha(\bar{t}), Y(\alpha(\bar{t}), \omega), a)$, Eq. \eqref{Pre-det-sys-for-BSDE(Mleshko)2 and 1} become
\begin{align}
      &\left(g\varphi_{y}-\varphi_s-\frac{1}{2}\varphi_{yy}Z^2\right)\left(\alpha(\bar{t}), Y(\alpha(\bar{t}), \omega), a\right)\label{Pre-det-sys-for-BSDE(Mleshko)before f version}\\
      &\nonumber =g\left(\bar{t}, \varphi\left(\alpha(\bar{t}), Y(\alpha(\bar{t}), \omega), a\right),\left(\frac{\varphi_y Z}{\eta}\right)\left(\alpha(\bar{t}), Y(\alpha(\bar{t}), \omega), a\right)\right)\eta^{2}\left(\alpha(\bar{t}), Y(\alpha(\bar{t}), \omega), a\right)
\end{align}
Substituting $\bar{t}=\beta(t)$ into Eq. \eqref{Pre-det-sys-for-BSDE(Mleshko)before f version}, one can write equation \eqref{Pre-det-sys-for-BSDE(Mleshko)before f version} as
\begin{align}
      &\nonumber\left(g\varphi_{y}-\varphi_t-\frac{1}{2}\varphi_{yy}Z^2\right)\left(t, Y(t, \omega), a\right)\\
      &=g\left(\beta(t), \varphi\left(t, Y(t, \omega), a\right),\left(\frac{\varphi_y Z}{\eta}\right)\left(t, Y(t, \omega), a\right)\right)\eta^{2}\left(t, Y(t, \omega), a\right)\label{Pre-det-sys-for-BSDE(Mleshko) finale version}
\end{align}
Differentiating Eq. \eqref{Pre-det-sys-for-BSDE(Mleshko) finale version} with respect to the parameter $a$, one obtains the following equation
\begin{align}
      &\left(g\varphi_{ay}-\varphi_{at}-\frac{1}{2}\varphi_{ayy}Z^2\right)\left(t, Y(t, \omega), a\right)\nonumber\\
      &=2\eta\eta_a\left(t, Y(t, \omega), a\right)g\left(\beta(t), \varphi\left(t, Y(t, \omega), a\right),\left(\frac{\varphi_y Z}{\eta}\right)\left(t, Y(t, \omega), a\right)\right)\nonumber\\
      &+\eta^2\left(t, Y(t, \omega), a\right)\left[g_t\frac{\partial\beta}{\partial a}+g_y\varphi_a+g_z\frac{\varphi_{ay}\eta-\eta_a\varphi_y}{\eta^2}Z\right]
      \label{det-sys-for-BSDE(Mleshko)before finale version}
\end{align}
Substituting $a=0$ into last Eq. \eqref{det-sys-for-BSDE(Mleshko)before finale version}, one has 
\begin{align}\label{Before det-sys-Meleshko}
 g\left(\left.\frac{\partial \varphi}{\partial a}\right|_{a=0}\right)_{y}-\left(\left.\frac{\partial \varphi}{\partial a}\right|_{a=0}\right)_{t}-\frac{1}{2} \left(\left.\frac{\partial \varphi}{\partial a}\right|_{a=0}\right)_{yy}Z^{2}&=2g\left.\frac{\partial\eta}{\partial a}\right|_{a=0}+g_t\left.\frac{\partial\beta}{\partial a}\right|_{a=0}+g_y\left.\frac{\partial\varphi}{\partial a}\right|_{a=0}\\ \nonumber
&+g_z\left(\left(\left.\frac{\partial\varphi}{\partial a}\right|_{a=0}\right)_y-\left.\frac{\partial\eta}{\partial a}\right|_{a=0}\right)Z
\end{align}
Since $\beta(t, \omega, a)=\int_{0}^{t} \eta^{2}(s, X(s, \omega), a) \mathrm{d} s$ for all $t \geqslant 0$, differentiating this with respect to $a$, one finds

\begin{equation}
\left.\frac{\partial \beta}{\partial a}\right|_{a=0}=\left.2 \int_{0}^{t} \frac{\partial \eta}{\partial a}\right|_{a=0} \mathrm{~d} s=2 \int_{0}^{t} \tau(s,y) \mathrm{~d} s\label{beta_a(a=0)=2inteta_a(a=0)=2int*tauds}
\end{equation}

Substituting $\left.\frac{\partial \beta}{\partial a}\right|_{a=0}$ into Eq \eqref{Before det-sys-Meleshko}, one arrives at the following equations

\begin{align}
    & \nonumber g\gamma_y(t,Y(t,\omega))-\gamma_t(t,Y(t,\omega))-\frac{1}{2}\gamma_{yy}(t,Y(t,\omega))Z^2(t,\omega)\\
    &\nonumber=2g\tau(t,Y(t,\omega))+ 2\left(\int_{0}^{t} \tau(s,Y(s,\omega)) \mathrm{~d} s\right) g_t
    +g_y\gamma(t,Y(t,\omega))+g_z\left(\gamma_y(t,Y(t,\omega)-\tau(t,Y(t,\omega))\right)Z\\
    \label{det-sys-meleshko final version}
\end{align}

Eq. (5.1.22) is an integro-differential equations for the functions $\tau (t, y)$ and $\gamma(t, y)$. These equa-
tion have to be satisfied for any solution $Y (t, \omega)$ of the backward stochastic differential equation \eqref{BSDE au sens de Meleshko}
\begin{defin}
A Lie group of transformations \eqref{lie group meleshko for Y(BSDE)} and \eqref{lie group meleshko for time t(BSDE)} is called admitted by the BSDE \eqref{BSDE au sens de Meleshko}, if $\gamma(t, x)$ and $\tau(t, x)$ satisfy the determining equation \eqref{det-sys-meleshko final version}.
\end{defin}

The determining equations for an admitted Lie group of transformations were derived under the assumption that the transformations \eqref{lie group meleshko for Y(BSDE)} and\eqref{lie group meleshko for time t(BSDE)} satisfy the requirements \eqref{det-sys-meleshko final version}.

Assume that one has found the functions $\tau(t, y)$ and $\gamma(t, y)$ which are solutions of the determining equations \eqref{det-sys-meleshko final version}. Then the Lie group of transformations \eqref{lie group meleshko for Y(BSDE)}, \eqref{lie group meleshko for time t(BSDE)} is recovered by solving the Lie equations
$$
\begin{aligned}
& \frac{\partial \Xi}{\partial a}(t, y, a)=h(\Xi, \varphi), \\
& \frac{\partial \varphi}{\partial a}(t, y, a)=\gamma(\Xi, \varphi),
\end{aligned}
$$
with the initial conditions,
$$
H(t, y, 0)=t, \quad \varphi(t, y, 0)=y
$$
where $h(t, y)=2 \int_{0}^{t} \tau(s, y) \mathrm{d} s$ and $
\eta^{2}=\frac{\partial \Xi}{\partial t}.
$
\subsection{Application to quadratic BSDE}
Inspired by Barrieu and El Karoui in \cite{10.1214/12-AOP743}, who studied the existence of solutions for BSDEs with a generator that is quadratic in $z$ (specifically, $g(z)=z^2)$, we consider the quadratic BSDE:
\begin{equation}
     Y(t,\omega) =Y(T,\omega)+\int_{t}^{T}Z^2_s ds-\int_{t}^{T}Z_s d B_{s}\label{quadratic BSDEs}
\end{equation}
Replacing $g$ in the determining equations \eqref{det-sys-meleshko final version}, we have: $$\left(-\frac{1}{2}\gamma_{yy}-\gamma_y\right)Z^2-\gamma_t=0.$$ Then, 
$$-\frac{1}{2}\gamma_{yy}-\gamma_y=0,\quad \gamma_t=0. $$
Solving yields,
$\gamma(t,y)=c_1\exp{(-2y)}+c_2$ and $\tau(t,y)=f(t,y)$, where $f(t,y)$ is an arbitrary function that is smooth enough and $c_1$ and $c_2$ are arbitrary constants.
\begin{remark}
\begin{itemize}
    \item 
          Since $\tau$ is an arbitrary function of $t$ and $y$, this suggests that for every choice of $\tau$, we obtain a new Brownian motion. This result is quite remarkable, though somewhat unexpected.
          \item 
          Zhang and Guangyan, in their work \cite{sym11091153}, computed the projectable symmetries of BSDEs using Gaeta/Kozlov approach in the quadratic case $g=C|z|^2$. They identified the following symmetries: $$\tau(t)=f(t),\quad \xi(t,x)=C_1e^{-2Cy}+C_2$$
          where $f$ is an arbitrary function smooth enough. This result aligns perfectly with what we found using the approach of \cite{srihirun2007definition}, which focuses on time transformations that depend on both time $t$ and the processes $Y$.
    \end{itemize}
\end{remark}
For example, consider $\tau(t,y)=f(t,y)=C_3$.\\
We have that $h(t,y)=2\int_0^t\tau(s,y)ds=2C_3t.$
For example, consider $C_1=C_3=1$ and $C_2=0$, we find the following symmetry $$v=\exp{(-2y)}\frac{\partial}{\partial y}+2t\frac{\partial}{\partial t}$$
We exponentiate the symmetry $v$ by solving Lie's equation \eqref{lie group meleshko for time t(BSDE)} and \eqref{lie group meleshko for Y(BSDE)}
\begin{equation}
    \frac{\partial \Xi}{\partial a}=2\Xi,\quad\quad \frac{\partial \varphi}{\partial a}=\exp{(-2\varphi)}
\end{equation}
with initial conditions 
\begin{equation}
    \Xi(t,x,0)=t,\quad\quad \varphi(t,x,0)=x 
\end{equation}

Solving the equations gives the set of transformations that acts on time $t$ and the process $y$: 
\begin{align}
    \bar{t}=\Xi(t,y,a)=e^{2a}t \label{bar(t) Meleshko BSDE example}\\
    \bar{y}=\varphi(t,y,a)=\frac{1}{2}\ln{(2a+e^{2y})}.\label{bar(y) Meleshko BSDE example}
\end{align}
Then $\eta^2=\Xi_t=e^{2a}$. Using \eqref{bar((Z) meleshko BSDE} we find the Lie group of transformation that acts on the process $z$ which is 
\begin{align}
    \bar{z}=\zeta(t,y,z,a)=e^{-a}\frac{e^{2y}}{2a+e^{2y}}z.\label{bar(z) Meleshko BSDE example}
\end{align}
Let us verify that the Lie group of transformation \eqref{bar(z) Meleshko BSDE example}, \eqref{bar(t) Meleshko BSDE example} and \eqref{bar(y) Meleshko BSDE example} transforms a solution of equation \eqref{BSDE au sens de Meleshko} into a solution of the same equation. Suppose that $(Y_t,Z_t)$ is a solution of equation \eqref{BSDE au sens de Meleshko}. Due to Theorem \ref{new bar(B) Meleshko}, the Brownian motion $B(t)$ is transformed to the Brownian motion 
\begin{align}
\Bar{B}_t=\int_0^te^adB_s=e^a B_t,\quad t\geq 0,\label{bar(B)-meleshko-exemple}
\end{align}
\begin{align}
    \text{where }\quad \beta(t)=\int_0^te^{2a}ds=te^{2a},\quad \alpha(t)=\inf_{s\geq 0}\{s:\beta(s)>t\}=\beta^{-1}(t),
\end{align}
According to Theorem \ref{new bar(B) Meleshko}, the process $\Bar{B}_{\alpha(t)}=e^aB_{\frac{t}{e^{2a}}}$ is a Brownian motion, which is a direct consequence of the scaling property.\\
Applying Ito formula to $\varphi(t,Y(t,\omega),a)=\frac{1}{2}\ln{(2a+e^{2Y(t,\omega)})}$, one finds that
\begin{align}
    \bar{Y}_T&=\frac{1}{2}\ln{\left(2a+e^{2Y(T,\omega)}\right)}=\frac{1}{2}\ln{\left(2a+e^{2Y(t,\omega)}\right)}+\int_t^T \frac{e^{2Y(t,\omega)}}{2a+e^{2Y(t,\omega)}}dY(s,\omega)\\
    &\nonumber+\int_t^T\frac{2a e^{2Y(t,\omega)}}{(2a+e^{2Y(t,\omega)})^2}d<Y>_s,\\
    &=\frac{1}{2}\ln{\left(2a+e^{2Y(t,\omega)}\right)}+\int_t^T \frac{e^{2Y(s,\omega)}}{2a+e^{2Y(s,\omega)}}d(-Z^2((s,\omega))ds+Z((s,\omega)dB(s))\\
    &\nonumber+\int_t^T\frac{2a e^{2Y(t,\omega)}}{(2a+e^{2Y(s,\omega)})^2}Z^2(s,\omega)ds,\\
    &=\frac{1}{2}\ln{\left(2a+e^{2Y(t,\omega)}\right)}+\int_t^T \left(\frac{-e^{2Y(s,\omega)}}{2a+e^{2Y(s,\omega)}}+\frac{2ae^{2Y(s,\omega)}}{(2a+e^{2Y(s,\omega)})^2}\right)Z^2((s,\omega)ds\label{Meleshko exmp}\\
    &\nonumber+\int_t^T\frac{e^{2Y(s,\omega)}}{2a+e^{2Y(s,\omega)}}Z(s,\omega)dB(s),
\end{align}
Using \eqref{alpha_bar(t)(s)}, $\frac{\partial \alpha}{\partial \bar{t}}=e^{2a}$. Changing the variable of the integral $s=\alpha(\bar{s})$ in the Reimann integral in \eqref{Meleshko exmp}, it becomes 
\begin{align}
    \int_t^T \left(\frac{-e^{2Y(s,\omega)}}{2a+e^{2Y(s,\omega)}}+\frac{2ae^{2Y(s,\omega)}}{((2a+e^{2Y(s,\omega))^2}}\right)Z^2((s,\omega)ds=\int_{\beta(t)}^{\beta(T)} \left(\frac{-e^{4Y(\alpha(\bar{s}),\omega)}e^{{-2a}}}{(2a+e^{2Y(\alpha(\bar{s}),\omega))^2}}\right)Z^2(\alpha(\bar{s}),\omega)d\alpha(\bar{s})
\end{align}
Because of the transformation of the Brownian motion \eqref{bar(B)-meleshko-exemple} and the change of Brownian motion in the Ito integral \eqref{time change formula oksendal}, the Ito integral in \eqref{Meleshko exmp} becomes
\begin{align}
    \int_t^T\frac{e^{2Y(s,\omega)}}{2a+e^{2Y(s,\omega)}}Z(s,\omega)dB(s)&=\int_{\beta(t)}^{\beta(T)}\frac{e^{-a}e^{2Y(\alpha(\bar{s}),\omega)}}{2a+e^{2Y(\alpha(\bar{s}),\omega)}}Z(\alpha(\bar{s}),\omega)d\bar{B}(\alpha(\bar{s}))
\end{align}
\begin{align}
   \text{Since}\quad \bar{Z}(\alpha(\bar{s}),\omega)&=\zeta(\alpha(\bar{s}),Y(\alpha(\bar{s}),\omega),Z(\alpha(\bar{s}),\omega),a))\\
   &=e^{-a}\frac{e^{2Y(\alpha(\bar{s}),\omega)}}{2a+e^{2Y(\alpha(\bar{s}),\omega)}}Z(\alpha(\bar{s}),\omega)
\end{align} 
\begin{align}
    \text{and}\quad Y(\bar{t},\omega)=\varphi(\alpha(\bar{t}),Y(\alpha(\bar{t}),\omega),a)=\frac{1}{2}\ln{\left(2a+e^{2Y(\alpha(\bar{t}),\omega)}\right)}
\end{align}
 one finds that 
\begin{align}
    Y(T,\omega)&=\frac{1}{2}\ln{\left(2a+e^{2Y(t,\omega)}\right)}-\int_{\beta(t)}^{\beta(T)} \left(\frac{e^{2Y(s,\omega)}e^{{-a}}}{2a+e^{2Y(s,\omega)}}Z(s,\omega)\right)^2ds\\
    &\nonumber+\int_{\beta(t)}^{\beta(T)}\left(\frac{e^{-a}e^{2Y(s,\omega)}}{2a+e^{2Y(s,\omega)}}Z(s,\omega)\right)d\bar{B}(s),\\
    &=\frac{1}{2}\ln{\left(2a+e^{2Y(t,\omega)}\right)}-\int_{\beta(t)}^{\beta(T)} \bar{Z}(s,\omega)^2ds+\int_{\beta(t)}^{\beta(T)}\bar{Z}(s,\omega)d\bar{B}(s)  
\end{align}
Because $\bar{Z}(\beta(t),\omega)=\zeta(t,Y(t,\omega),Z(t,\omega),a),\,\,\text{and}\,\, \bar{Y}(\beta(t),\omega)=\varphi(t,Y(t,\omega),a)$
one has,
\begin{align}
     &\bar{Y}(\beta(t),\omega)=\frac{1}{2}\ln{\left(2a+e^{2Y(\beta(t),\omega)}\right)}=\frac{1}{2}\ln{\left(2a+e^{2Y(T,\omega)}\right)}+\int_{\beta(t)}^{\beta(T)} \bar{Z}(s,\omega)^2ds-\int_{\beta(t)}^{\beta(T)}\bar{Z}(s,\omega)d\bar{B}(s)  
\end{align}
This confirms that the Lie group of transformations \eqref{lie group meleshko for Y(BSDE)}, \eqref{lie group meleshko for time t(BSDE)} and \eqref{bar((Z) meleshko BSDE} transforms any solution of Eq. \eqref{BSDE au sens de Meleshko} into a solution of the same equation, but with a different terminal condition.
\begin{remark}
    $(\bar{Y},\bar{Z})$ is a solution of the quadratic BSDE \eqref{quadratic BSDEs} driven by the Brownian motion $\bar{B}_t$ but they are not adapted to the filtration of $(\bar{B}_t)_{_{0\leq t\leq T }}$. They are adapted to $\mathcal{F}_{\alpha(t)}$ which contains $\sigma(\bar{B}_s,\, s\leq t)$. Moreover, using the martingale representation theorem \cite{revuz2013continuous}, any $\mathcal{F}_{\alpha(t)}$ square integrable martingale $M_t$ can be represented as a stochastic integral with respect to $\bar{B}_t$ $$M_t=\mathbb{E}(M_0)+\int_0^t \psi(s)d\bar{B}_s$$ where $\psi(.)$ is $\mathcal{F}_{\alpha(t)}$ predictable.
\end{remark}

\section{Admitted symmetries of forward backward stochastic differential equations}
In this section, we will compute the admitted symmetries of the following  uncoupled FBSDE:
\begin{align}
    \begin{cases}
       X(t,\omega) &=X(0,\omega)+ \displaystyle \int_0^t b(s, X(s,\omega)) ds +\int_0^t \sigma(s, X(s,\omega)) dW_s, \\
        Y(t,\omega)&=Y(T,\omega)+\displaystyle \int_t^T g\left(s, X(s,\omega), Y(s,\omega), Z(s,\omega)\right) ds -\int_t^T Z(s,\omega) dW_s, \\\label{FBSDE int form Meleshko}
         Y(T,\omega) &= H(X(T,\omega)).
    \end{cases}
\end{align}
Here, the Lie groups of transformations that acts on the process $X$, $Y$ and time $t$ are the following: 
\begin{align}
 \bar{x}&=\phi(t,x,y,a)\label{lie group meleshko for X(FBSDE)}\\
    \bar{y}&=\varphi(t,x,y,a)\label{lie group meleshko for Y(FBSDE)}\\
    \bar{t}&=\Xi(t,x,y,a)\label{lie group meleshko for time t(FBSDE)}
\end{align}
As seen in the admitted symmetries of BSDE, the Lie group that transforms $Z$ will be derived from the transformations of $Y$ and time $t$. The infinitesimal generator that generates an admitted symmetry of \eqref{FBSDE int form Meleshko} is of the form
\begin{equation}
    h(t,x,y)\partial_t+\xi(t,x,y)\partial_x+\gamma(t,x,y)\partial_y
\end{equation}
such that
 \begin{align}
      h(t,x,y)&=\frac{\partial \Xi(t,x,y,a)}{\partial a}|_{a=0}=\Xi_a(t,x,y,0),\\
     \xi(t,x,y)&=\frac{\partial \phi(t,x,y,a)}{\partial a}|_{a=0}=\phi_a(t,x,y,0),\\
     \gamma(t,x,y)&=\frac{\partial \varphi(t,x,y,a)}{\partial a}|_{a=0}=\varphi_a(t,x,y,0),
 \end{align}
and $\Xi$, $\phi$ and $\varphi$ satisfy Lie's equations 
\begin{equation}
    \frac{\partial \Xi}{\partial a}=h(\Xi,\phi,\varphi),\quad  \frac{\partial \varphi}{\partial a}=\gamma(\Xi,\phi,\varphi),\quad \frac{\partial \phi}{\partial a}=\xi(\Xi,\phi,\varphi)
\end{equation}
with initial conditions for $a=0$,
\begin{equation}
    \Xi(t,x,y,0)=t,\quad \varphi(t,x,y,0)=y,\quad \phi(t,x,y,0)=x.\label{initial cond for lie equ FBSDE}
\end{equation}

We define 
\begin{align}
\Bar{X}(\Bar{t},\omega)&=\phi(\alpha(\Bar{t}),X(\alpha(\Bar{t}),\omega),Y(\alpha(\Bar{t}),\omega),a),\\
\Bar{Y}(\Bar{t},\omega)&=\varphi(\alpha(\Bar{t}),X(\alpha(\Bar{t}),\omega),Y(\alpha(\Bar{t}),\omega),a),
\end{align}
 the transformations of the process $X$ and $Y$. 
 For the transformed Brownian motion, we again apply Theorem \eqref{new bar(B) Meleshko}; however, in this case, the time change rate 
$\eta$ also depends on $X_t$, the solution of the forward equation in \eqref{FBSDE int form Meleshko}.
\begin{thm}
    Let $\eta(t,x,y,a)$ and $Y_t$ and $X_t$ be as above and $B_t$ a standard Brownian motion. Define 
    $$\Bar{B}_t=\int_0^t\eta(s,Y(s,\omega),X(s,\omega),a)dB_s,\quad\quad t\geq 0$$
    Then $(\Bar{B}_{\alpha(t)},\mathcal{F}_{\alpha(t)})$ is a standard Brownian motion, where \\
$\mathcal{F}_{\alpha(t)}=\{A\in\mathcal{F}:A\cap\{\omega:\alpha(t)\leq s\}\in\mathcal{F}_s\quad \text{for all}\,\, s\geq 0\}$ and $\alpha$ as defined above.
\end{thm}
\noindent For a complete proof of the theorem, 
 we refer to \cite{srihirun2007definition}. Unlike Theorem \ref{new bar(B) Meleshko}, here the rate $\eta$ also depends on the process $X$ in \ref{FBSDE int form Meleshko}.
 
 \subsection{Determining equations of FBSDE}
Suppose that $(X(t, \omega),Y(t, \omega),Z(t, \omega))$ satisfy the FBSDE \eqref{FBSDE int form Meleshko}. \\
Define, $$\chi(t,\omega)=\Bar{Y}(\beta(t),\omega) \,\,\text{and}\,\, \psi(t,\omega)=\Bar{X}(\beta(t),\omega)$$\\
For almost all $\omega$,
\begin{align*}
\psi(t, \omega)&=\phi \left(\alpha(\beta(t)),X(\alpha(\beta(t),\omega),Y(\alpha(\beta(t)), \omega), a\right)=\phi(t,X(t, \omega), Y(t, \omega), a).\\
    \chi(t, \omega)&=\varphi \left(\alpha(\beta(t)),X(\alpha(\beta(t),\omega),Y(\alpha(\beta(t)), \omega), a\right)=\varphi(t,X(t, \omega), Y(t, \omega), a).
\end{align*}
By applying Itô’s formula, we have   
\begin{align}
\nonumber \psi(t, \omega)&=\phi(0,X(0, \omega), Y(0, \omega), a)+\int_{0}^{t}
\phi_{x}(s,X(s,\omega),Y(s,\omega),a)dX(s,\omega)\\
&+\int_{0}^{t} \phi_{y}(s,Y(s,\omega),a)dY(s,\omega),a)
+\int_{0}^{t} \phi_{s}(s,X(s,\omega),Y(s,\omega),a)ds\\
&\nonumber+\frac{1}{2}\Bigg[\int_{0}^{t}
\phi_{xx}(s,X(s,\omega),Y(s,\omega),a)d<X,X>_s+\phi_{yy}(s,X(s,\omega),Y(s,\omega),a)d<Y,Y>_s\\
&\nonumber +\phi_{xy}(s,X(s,\omega),Y(s,\omega),a)d<X,Y>_s\Bigg]\\
&=\nonumber\phi(0,X(0, \omega), Y(0, \omega), a)+\\
&\nonumber\int_{0}^{t}
\phi_{x}(s,X(s,\omega),Y(s,\omega),a)[b(s,X(s,\omega)) d s+\sigma(s,X(s,\omega)) d B(s)]\\
&+\int_{0}^{t} \phi_{y}(s,Y(s,\omega),a)\left[-g(s,X(s,\omega),Y(s,\omega),Z(s,\omega)) d s+Z(s,\omega)d B(s)\right]\\
&\nonumber+\int_{0}^{t} \phi_{s}(s,X(s,\omega),Y(s,\omega),a)ds+\frac{1}{2}\Bigg[\int_{0}^{t}
\phi_{xx}(s,X(s,\omega),Y(s,\omega),a)\sigma^2(s,X(s,\omega))ds\\
&\nonumber+\phi_{yy}(s,X(s,\omega),Y(s,\omega),a)Z^2(s,\omega)ds
+\phi_{xy}(s,X(s,\omega),Y(s,\omega),a)\sigma(s,X(s,\omega))Z(s,\omega)ds\Bigg]\\
&=\phi(0,X(0, \omega), Y(0, \omega), a)+\int_{0}^{t} [\sigma\phi_x+Z\phi_y](s,X(s,\omega),Y(s,\omega))dB(s)\\
&+ \nonumber\int_{0}^{t}
[b\phi_{x}-g\phi_y+\phi_s+\frac{1}{2}
\sigma^2\phi_{xx}+\frac{1}{2}
Z^2\phi_{yy}+
\sigma Z\phi_{xy}](s,X(s,\omega),Y(s,\omega))ds.
\end{align}
 and 
 \begin{align}
\nonumber \chi(t, \omega)&=\varphi(T,X(T, \omega), Y(T, \omega), a)-\int_{t}^{T}
\varphi_{x}(s,X(s,\omega),Y(s,\omega),a)dX(s,\omega)\\
&-\int_{t}^{T} \varphi_{y}(s,Y(s,\omega),a)dY(s,\omega),a)
-\int_{t}^{T} \varphi_{s}(s,X(s,\omega),Y(s,\omega),a)ds\\
&\nonumber-\frac{1}{2}\Bigg[\int_{t}^{T}
\varphi_{xx}(s,X(s,\omega),Y(s,\omega),a)d<X,X>_s+\varphi_{yy}(s,X(s,\omega),Y(s,\omega),a)d<Y,Y>_s\\
&\nonumber +\varphi_{xy}(s,X(s,\omega),Y(s,\omega),a)d<X,Y>_s\Bigg]\\
&=\varphi(T,X(T, \omega), Y(T, \omega), a)-\int_{t}^{T} [\sigma\varphi_x+Z\varphi_y](s,X(s,\omega),Y(s,\omega),a)dB(s)\\
&\nonumber-\int_{t}^{T}
[b\varphi_{x}-g\varphi_y+\varphi_s+\frac{1}{2}
\sigma^2\varphi_{xx}+\frac{1}{2}
Z^2\varphi_{yy}+
\sigma Z\varphi_{xy}](s,X(s,\omega),Y(s,\omega),a)ds.
\end{align}
 By applying the random time change formula as established in the previous paragraph \ref{Adm sym of BSDE (paragraph)}, we obtain 
 \begin{align}
     \psi(t,\omega)&=\phi(0,X(0, \omega), Y(0, \omega), a)+\int_{0}^{\beta(t)} \frac{[\sigma\phi_x+Z\phi_y]}{\eta}(\alpha(s),X(\alpha(s),\omega),Y(\alpha(s),\omega),a)d\bar{B}(s)\label{psi FBSDE Meleshko oks change of var}\\
+ \nonumber\int_{0}^{\beta(t)}&
\left[b\phi_{x}-g\phi_y+\phi_s+\frac{1}{2}
\sigma^2\phi_{xx}+\frac{1}{2}
Z^2\phi_{yy}+
\sigma Z\phi_{xy}\right](\alpha(s),X(\alpha(s),\omega),Y(\alpha(s),\omega),a)\alpha_{\bar{t}}(s)ds.
 \end{align}
 and 
 \begin{align}
     \chi(t,\omega)&=\varphi(T,X(T, \omega), Y(T, \omega), a)-\int_{\beta(t)}^{\beta(T)} \frac{[\sigma\varphi_x+Z\varphi_y]}{\eta}(\alpha(s),X(\alpha(s),\omega),Y(\alpha(s),\omega),a)d\bar{B}(s)\label{chi FBSDE Meleshko oks change of var}\\
- \nonumber\int_{\beta(t)}^{\beta(T)}&
\left[b\varphi_{x}-g\varphi_y+\varphi_s+\frac{1}{2}
\sigma^2\varphi_{xx}+\frac{1}{2}
Z^2\varphi_{yy}+
\sigma Z\varphi_{xy}\right](\alpha(s),X(\alpha(s),\omega),Y(\alpha(s),\omega),a)\alpha_{\bar{t}}(s)ds.
 \end{align}
 Substituting \eqref{alpha_bar(t)(s)} into equations \eqref{psi FBSDE Meleshko oks change of var} and \eqref{chi FBSDE Meleshko oks change of var} 
  \begin{align}
     \psi(t,\omega)&=\phi(0,X(0, \omega), Y(0, \omega), a)+\int_{0}^{\beta(t)} \frac{[\sigma\phi_x+Z\phi_y]}{\eta}(\alpha(s),X(\alpha(s),\omega),Y(\alpha(s),\omega),a)d\bar{B}(s)\\
+ \nonumber\int_{0}^{\beta(t)}&
\frac{\left[b\phi_{x}-g\phi_y+\phi_s+\frac{1}{2}
\sigma^2\phi_{xx}+\frac{1}{2}
Z^2\phi_{yy}+
\sigma Z\phi_{xy}\right]}{\eta^2}(\alpha(s),X(\alpha(s),\omega),Y(\alpha(s),\omega),a)ds.
 \end{align}
 and 
 \begin{align}
     \chi(t,\omega)&=\varphi(T,X(T, \omega), Y(T, \omega), a)-\int_{\beta(t)}^{\beta(T)} \frac{[\sigma\varphi_x+Z\varphi_y]}{\eta}(\alpha(s),X(\alpha(s),\omega),Y(\alpha(s),\omega),a)d\bar{B}(s)\\
- \nonumber\int_{\beta(t)}^{\beta(T)}&
\frac{\left[b\varphi_{x}-g\varphi_y+\varphi_s+\frac{1}{2}
\sigma^2\varphi_{xx}+\frac{1}{2}
Z^2\varphi_{yy}+
\sigma Z\varphi_{xy}\right]}{\eta^2}(\alpha(s),X(\alpha(s),\omega),Y(\alpha(s),\omega),a)ds.
 \end{align}
 Requiring that transformations \eqref{lie group meleshko for time t(FBSDE)}, \eqref{lie group meleshko for X(FBSDE)} and \eqref{lie group meleshko for Y(FBSDE)} map a solution of the FBSDE \eqref{FBSDE int form Meleshko} into a solution of the same FBSDE, one obtains
\begin{align}
    \begin{cases}
       \bar{X}(\bar{t},\omega) &=\bar{X}(0,\omega)+ \displaystyle \int_0^{\bar{t}} b(s, \bar{X}(s,\omega)) ds +\int_0^{\bar{t}} \sigma(s, \bar{X}(s,\omega)) d\bar{B}_s, \\
        \bar{Y}(\bar{t},\omega)&=\bar{Y}(\bar{T},\omega)+\displaystyle \int_{\bar{t}}^{\bar{T}} g\left(s, \bar{X}(s,\omega), \bar{Y}(s,\omega), \bar{Z}(s,\omega)\right) ds -\int_{\bar{t}}^{\bar{T}} \bar{Z}(s,\omega) d\bar{B}_s, \\
         \bar{Y}(\bar{T},\omega) &= H(\bar{X}(\bar{T},\omega)).
    \end{cases}
\end{align}
 Replacing $\bar{t}=\beta(t)$ and $\bar{T}=\beta(T)$ into these equation, one gets
 \begin{align}
    \begin{cases}
       \bar{X}(\beta(t),\omega) &=\bar{X}(0,\omega)+ \displaystyle \int_0^{\beta(t)} b(s, \bar{X}(s,\omega)) ds +\int_0^{\beta(t)} \sigma(s, \bar{X}(s,\omega)) d\bar{B}_s, \\
        \bar{Y}(\beta(t),\omega)&=\bar{Y}(\beta(T),\omega)+\displaystyle \int_{\beta(t)}^{\beta(T)} g\left(s, \bar{X}(s,\omega), \bar{Y}(s,\omega), \bar{Z}(s,\omega)\right) ds -\int_{\beta(t)}^{\beta(T)} \bar{Z}(s,\omega) d\bar{B}_s, \\
         \bar{Y}(\beta(T),\omega) &= H(\bar{X}(\beta(T),\omega)).
    \end{cases}
\end{align}
Comparing the Riemann integrals and the Itô integrals, we obtain
\begin{align}
    &b(s, \bar{X}(s,\omega))=\frac{\left[b\phi_{x}-g\phi_y+\phi_s+\frac{1}{2}
\sigma^2\phi_{xx}+\frac{1}{2}
Z^2\phi_{yy}+
\sigma Z\phi_{xy}\right]}{\eta^2}(\alpha(s),X(\alpha(s),\omega),Y(\alpha(s),\omega),a)\label{b(s, bar(X)(s,omega))=...FBSDE-Meleshko}\\
&\sigma(s, \bar{X}(s,\omega))=\frac{[\sigma\phi_x+Z\phi_y]}{\eta}(\alpha(s),X(\alpha(s),\omega),Y(\alpha(s),\omega),a)\label{sigma(s, bar(X)(s,omega))=...FBSDE-Meleshko}\\
&g\left(s, \bar{X}(s,\omega), \bar{Y}(s,\omega), \bar{Z}(s,\omega)\right)=\label{g(s,bar(X),bar(Y),bar(Z))-FBSDE-Meleshko}\\
&\nonumber-\frac{\left[b\varphi_{x}-g\varphi_y+\varphi_s+\frac{1}{2}
\sigma^2\varphi_{xx}+\frac{1}{2}
Z^2\varphi_{yy}+
\sigma Z\varphi_{xy}\right]}{\eta^2}(\alpha(s),X(\alpha(s),\omega),Y(\alpha(s),\omega),a)\\
&\bar{Z}(s,\omega)=\frac{[\sigma\varphi_x+Z\varphi_y]}{\eta}(\alpha(s),X(\alpha(s),\omega),Y(\alpha(s),\omega),a)\label{bar(Z)=(sigma*varphi_x+Z*varphi_y)/eta-FBSDE-Meleshko}
 \end{align}
We can now give the Lie group that transforms the variable $Z$, which we shall refer to as $\zeta$:
\begin{prop}
     \begin{equation}
     \zeta:=\bar{Z}=\frac{\varphi_y Z+\sigma\varphi_x}{\eta}.
 \end{equation}
\end{prop}

 Replacing $\eqref{bar(Z)=(sigma*varphi_x+Z*varphi_y)/eta-FBSDE-Meleshko}$ in $\eqref{g(s,bar(X),bar(Y),bar(Z))-FBSDE-Meleshko}$, we get 
 \begin{align}
     &g\left(s, \bar{X}(s,\omega), \bar{Y}(s,\omega),\frac{[\sigma\varphi_x+Z\varphi_y]}{\eta}(\alpha(s),X(\alpha(s),\omega),Y(\alpha(s),\omega),a)\right)=\label{g(s,bar(X),bar(Y),Z)-FBSDE-Meleshko}\\
&\nonumber-\frac{\left[b\varphi_{x}-g\varphi_y+\varphi_s+\frac{1}{2}
\sigma^2\varphi_{xx}+\frac{1}{2}
Z^2\varphi_{yy}+ 
\sigma Z\varphi_{xy}\right]}{\eta^2}(\alpha(s),X(\alpha(s),\omega),Y(\alpha(s),\omega),a)
 \end{align}
 Since $\bar{X}(\bar{t}, \omega)=\phi(\alpha(\bar{t}), Y(\alpha(\bar{t}),X(\alpha(\bar{t}) \omega), a)$ and $\bar{Y}(\bar{t}, \omega)=\varphi(\alpha(\bar{t}), Y(\alpha(\bar{t}),X(\alpha(\bar{t}) \omega), a)$, Eqs. \eqref{b(s, bar(X)(s,omega))=...FBSDE-Meleshko}, \eqref{sigma(s, bar(X)(s,omega))=...FBSDE-Meleshko} and \eqref{g(s,bar(X),bar(Y),Z)-FBSDE-Meleshko} become
 \begin{align}
    &b(\bar{t}, \phi(\alpha(\bar{t}), Y(\alpha(\bar{t}),X(\alpha(\bar{t}) \omega), a))\eta^2(\alpha(\bar{t}),X(\alpha(\bar{t}),\omega),Y(\alpha(\bar{t}),\omega),a)\label{pre-det-sysb(..)*eta^2=(...)}\\
    &\nonumber =\left[b\phi_{x}-g\phi_y+\phi_t+\frac{1}{2}
\sigma^2\phi_{xx}+\frac{1}{2}
Z^2\phi_{yy}+
\sigma Z\phi_{xy}\right](\alpha(\bar{t}),X(\alpha(\bar{t}),\omega),Y(\alpha(\bar{t}),\omega),a)\\
&\sigma(\bar{t}, \phi(\alpha(\bar{t}), Y(\alpha(\bar{t}),X(\alpha(\bar{t}), \omega), a))\eta(\alpha(\bar{t}),X(\alpha(\bar{t}),\omega),Y(\alpha(\bar{t}),\omega),a)\label{pre-det-sys-sigma*eta=sigma*phi_x+z*phi_y(FBSDE-Meleshko)}\\
&\nonumber=[\sigma\phi_x+Z\phi_y](\alpha(\bar{t}),X(\alpha(\bar{t}),\omega),Y(\alpha(\bar{t}),\omega),a)\\
&\nonumber g\Bigg(\bar{t}, \phi(\alpha(\bar{t}), Y(\alpha(\bar{t}),X(\alpha(\bar{t}) \omega), a), \varphi(\alpha(\bar{t}), Y(\alpha(\bar{t}),X(\alpha(\bar{t}) \omega), a),\\
&\frac{[\sigma\varphi_x+Z\varphi_y]}{\eta}(\alpha(\bar{t}),X(\alpha(\bar{t}),\omega),Y(\alpha(\bar{t}),\omega),a)\Bigg)\eta^2(\alpha(\bar{t}),X(\alpha(\bar{t}),\omega),Y(\alpha(\bar{t}),\omega),a)\label{pre-det-sys-g(..)*eta^2=-(...)(FBSDE-Meleshko)}\\
&\nonumber=-\left[b\varphi_{x}-g\varphi_y+\varphi_t+\frac{1}{2}
\sigma^2\varphi_{xx}+\frac{1}{2}
Z^2\varphi_{yy}+ 
\sigma Z\varphi_{xy}\right](\alpha(\bar{t}),X(\alpha(\bar{t}),\omega),Y(\alpha(\bar{t}),\omega),a)
 \end{align}
Substituting $\bar{t}=\beta(t)$ into Eqs.\eqref{pre-det-sys-sigma*eta=sigma*phi_x+z*phi_y(FBSDE-Meleshko)}, \eqref{pre-det-sys-g(..)*eta^2=-(...)(FBSDE-Meleshko)}, \eqref{pre-det-sysb(..)*eta^2=(...)} we can rewrite equations as 
\begin{align}
    &b(\beta(t), \phi(t, Y(t,\omega),X(t, \omega), a))\eta^2(t,X(t,\omega),Y(t,\omega),a))\\
    &\nonumber =\left[b\phi_{x}-g\phi_y+\phi_t+\frac{1}{2}
\sigma^2\phi_{xx}+\frac{1}{2}
Z^2\phi_{yy}+
\sigma Z\phi_{xy}\right](t,X(t,\omega),Y(t,\omega),a)\\
&\sigma(\beta(t), \phi(t, Y(t,\omega),X(t, \omega), a))\eta(t,X(t,\omega),Y(t,\omega),a)\\
&\nonumber=[\sigma\phi_x+Z\phi_y](t,X(t,\omega),Y(t,\omega),a)\\
&\nonumber g\Bigg(\beta(t), \phi(t, Y(t,\omega),X(t, \omega), a), \varphi(t, Y(t,\omega),X(t, \omega), a),\\
&\frac{[\sigma\varphi_x+Z\varphi_y]}{\eta}(t,X(t,\omega),Y(t,\omega),a)\Bigg)\eta^2(t,X(t,\omega),Y(t,\omega),a)\\
&\nonumber=-\left[b\varphi_{x}-g\varphi_y+\varphi_t+\frac{1}{2}
\sigma^2\varphi_{xx}+\frac{1}{2}
Z^2\varphi_{yy}+ 
\sigma Z\varphi_{xy}\right](t,X(t,\omega),Y(t,\omega),a)
 \end{align}
 Differentiating last equations with respect to the parameter $a$, one obtains the following equations
 \begin{align}
    &\left(2b\eta\eta_a+\eta^2\left(b_t\frac{\partial\beta}{\partial a}+b_x\phi_a\right)\right)(t,X(t,\omega),Y(t,\omega),a))\\
    &\nonumber =\left[b\phi_{xa}-g\phi_{ya}+\phi_{ta}+\frac{1}{2}
\sigma^2\phi_{xxa}+\frac{1}{2}
Z^2\phi_{yya}+
\sigma Z\phi_{xya}\right](t,X(t,\omega),Y(t,\omega),a)\\
&\left(\eta_a\sigma+\eta\left(\sigma_t\frac{\partial\beta}{\partial a}+\sigma_x\phi_a\right)\right)(t,X(t,\omega),Y(t,\omega),a)\\
&\nonumber=[\sigma\phi_{xa}+Z\phi_{ya}](t,X(t,\omega),Y(t,\omega),a)\\
&\nonumber\eta^2\left[g_t \frac{\partial\beta}{\partial a}+g_x\phi_a+g_y\varphi_a+g_z\left(\frac{\eta_a(\sigma\varphi_x+Z\varphi_y)-\eta(\sigma\varphi_{xa}+Z\varphi_{ya})}{\eta^2}\right)\right](t,X(t,\omega),Y(t,\omega),Z(t,\omega))\\
& +2\eta\eta_a g(t,X(t,\omega),Y(t,\omega),Z(t,\omega))\\
&\nonumber=-\left[b\varphi_{xa}-g\varphi_{ya}+\varphi_{ta}+\frac{1}{2}
\sigma^2\varphi_{xxa}+\frac{1}{2}
Z^2\varphi_{yya}+ 
\sigma Z\varphi_{xya}\right](t,X(t,\omega),Y(t,\omega), Z(t,\omega),a)
 \end{align}
 Substituting $a=0$, one has
  \begin{align}
    &\left( 2b\left.\frac{\partial\eta}{\partial a}\right|_{a=0}+b_t\left(\left.\frac{\partial \beta}{\partial a}\right|_{a=0}\right)+b_x\left(\left.\frac{\partial \phi}{\partial a}\right|_{a=0}\right)\right)=b\left(\left.\frac{\partial \phi}{\partial a}\right|_{a=0}\right)_x-g\left(\left.\frac{\partial \phi}{\partial a}\right|_{a=0}\right)_y\\
    &\nonumber+\left(\left.\frac{\partial \phi}{\partial a}\right|_{a=0}\right)_t+\frac{1}{2}
\sigma^2\left(\left.\frac{\partial \phi}{\partial a}\right|_{a=0}\right)_{xx}+\frac{1}{2}
Z^2\left(\left.\frac{\partial \phi}{\partial a}\right|_{a=0}\right)_{yy}+
\sigma Z\left(\left.\frac{\partial \phi}{\partial a}\right|_{a=0}\right)_{xy}\\
\nonumber\\
&\sigma \left.\frac{\partial\eta}{\partial a}\right|_{a=0}+\sigma_t\left.\frac{\partial\beta}{\partial a}\right|_{a=0}+\sigma_x\left.\frac{\partial\phi}{\partial a}\right|_{a=0}=\sigma\left(\left.\frac{\partial \phi}{\partial a}\right|_{a=0}\right)_x+Z\left(\left.\frac{\partial \phi}{\partial a}\right|_{a=0}\right)_y\\
\nonumber\\
&\nonumber g_t \left(\left.\frac{\partial \beta}{\partial a}\right|_{a=0}\right)+g_x\left(\left.\frac{\partial \phi}{\partial a}\right|_{a=0}\right)+g_y\left(\left.\frac{\partial \varphi}{\partial a}\right|_{a=0}\right)+g_z\left[\left(\left.\frac{\partial \eta}{\partial a}\right|_{a=0}-\left(\left.\frac{\partial \varphi}{\partial a}\right|_{a=0}\right)_y\right)Z-\sigma\left(\left.\frac{\partial \varphi}{\partial a}\right|_{a=0}\right)_x\right]\\
& +2\left(\left.\frac{\partial \eta}{\partial a}\right|_{a=0}\right) g=-\Bigg[b\left(\left.\frac{\partial \varphi}{\partial a}\right|_{a=0}\right)_x-g\left(\left.\frac{\partial \varphi}{\partial a}\right|_{a=0}\right)_y+\left(\left.\frac{\partial \varphi}{\partial a}\right|_{a=0}\right)_t\\
&\nonumber+\frac{1}{2}
\sigma^2\left(\left.\frac{\partial \varphi}{\partial a}\right|_{a=0}\right)_{xx}+\frac{1}{2}
Z^2\left(\left.\frac{\partial \varphi}{\partial a}\right|_{a=0}\right)_{yy}+ 
\sigma Z\left(\left.\frac{\partial \varphi}{\partial a}\right|_{a=0}\right)_{xy}\Bigg]
 \end{align}

\noindent Substituting $\left.\frac{\partial \beta}{\partial a}\right|_{a=0}$  expressed in \eqref{beta_a(a=0)=2inteta_a(a=0)=2int*tauds} into last equations, one arrives at the following equations:
\begin{align}\label{det-sys-final version-Meleshko-2b.tau}
    &  2b\tau(t,X(t,\omega),Y(t,\omega))+b_t 2 \int_{0}^{t} \tau(s,X(t,\omega),Y(t,\omega)) \mathrm{~d} s+b_x\xi(t,X(t,\omega),Y(t,\omega))\nonumber\\
    &+\xi_t(t,X(t,\omega),Y(t,\omega))
    =b\xi_x(t,X(t,\omega),Y(t,\omega))-g\xi_y(t,X(t,\omega),Y(t,\omega))\nonumber\\
    &+\frac{1}{2}
\sigma^2\xi_{xx}(t,X(t,\omega),Y(t,\omega))+\frac{1}{2}
Z^2\xi_{yy}(t,X(t,\omega),Y(t,\omega))+
\sigma Z\xi_{xy}(t,X(t,\omega),Y(t,\omega))
\end{align}
\begin{align}\label{det-sys-final version-Meleshko-sigma.tau}
&\sigma \tau(t,X(t,\omega),Y(t,\omega))+\sigma_t \int_{0}^{t} \tau(s,X(t,\omega),Y(t,\omega)) \mathrm{~d} s+\sigma_x \xi(s,X(t,\omega),Y(t,\omega))\nonumber\\
&=\sigma \xi_x(t,X(t,\omega),Y(t,\omega)) +Z\xi_y(t,X(t,\omega),Y(t,\omega))
\end{align}
\begin{align}\label{det-sys-final version-Meleshko-g_t2}
& g_t 2 \int_{0}^{t} \tau(s,X(t,\omega),Y(t,\omega)) \mathrm{~d} s+g_x  \xi(s,X(t,\omega),Y(t,\omega))+g_y \gamma(s,X(t,\omega),Y(t,\omega))\nonumber\\
&\nonumber+g_z \left(\tau(t,X(t,\omega),Y(t,\omega))-\gamma_y(t,X(t,\omega),Y(t,\omega))\right)Z-\sigma \gamma_x(t,X(t,\omega),Y(t,\omega))\nonumber\\
& +2\tau(t,X(t,\omega),Y(t,\omega)) g=-\Bigg[b\gamma_x(t,X(t,\omega),Y(t,\omega))-g\gamma_y(t,X(t,\omega),Y(t,\omega)\nonumber\\
&+\gamma_t(t,X(t,\omega),Y(t,\omega))
+\frac{1}{2}
\sigma^2\gamma_{xx}(t,X(t,\omega),Y(t,\omega))+\frac{1}{2}
Z^2\gamma_{yy}(t,X(t,\omega),Y(t,\omega))\nonumber\\
 &+\sigma Z\gamma_{xy}(t,X(t,\omega),Y(t,\omega))\Bigg]
 \end{align}
\begin{thm}\cite[p.~12]{sym11091153}
     For any $(X_t,Y_t,Z_t)$ solution of the FBSDE \eqref{FBSDE int form Meleshko}, we have that 
     \begin{equation}
         \Tilde{Y}_{\Tilde{T}} = H(\Tilde{X}_{\tilde{T}}) \,\,\text{holds if and olny if}\,\, \gamma(t,x,H(x))=H_x \xi(t,x,H(x))\label{terminal cond meleshko hold}
     \end{equation}
\end{thm}

Equations \eqref{det-sys-final version-Meleshko-2b.tau}, \eqref{det-sys-final version-Meleshko-sigma.tau}, \eqref{det-sys-final version-Meleshko-g_t2} and \eqref{terminal cond meleshko hold} are integro-differential equations for the functions $\gamma$, $\xi$ and $\tau$. These equations have to be satisfied for any solution $(X(t,\omega),Y(t,\omega),Z(t,\omega))$ of the forward-backward stochastic differential equation \eqref{FBSDE int form Meleshko}.

\begin{defin}
A Lie group of transformations \eqref{lie group meleshko for time t(FBSDE)}, \eqref{lie group meleshko for X(FBSDE)} and \eqref{lie group meleshko for Y(FBSDE)} are called admitted by the FBSDE \eqref{FBSDE int form Meleshko}, if $\gamma(t, x,y)$, $\tau(t, x,y)$ and $\xi(t,x,y)$ satisfy the determining equation  \eqref{det-sys-final version-Meleshko-2b.tau}, \eqref{det-sys-final version-Meleshko-g_t2}, \eqref{det-sys-final version-Meleshko-sigma.tau} and \eqref{terminal cond meleshko hold}.
\end{defin}

The determining equations for an admitted Lie group of transformations were derived under the assumption that the transformations \eqref{lie group meleshko for time t(FBSDE)}, \eqref{lie group meleshko for X(FBSDE)} and \eqref{lie group meleshko for Y(FBSDE)} satisfy the requirements \eqref{det-sys-final version-Meleshko-2b.tau}, \eqref{det-sys-final version-Meleshko-g_t2}, \eqref{det-sys-final version-Meleshko-sigma.tau} and \eqref{terminal cond meleshko hold}.

Assume that one has found the functions $\tau(t, x,y)$, $\gamma(t,x,y)$ and $\xi(t,x,y)$ which are solutions of the determining equations \eqref{det-sys-final version-Meleshko-2b.tau}, \eqref{det-sys-final version-Meleshko-g_t2}, \eqref{det-sys-final version-Meleshko-sigma.tau} and \eqref{terminal cond meleshko hold}. Then the Lie group of transformations \eqref{lie group meleshko for time t(FBSDE)}, \eqref{lie group meleshko for X(FBSDE)} and \eqref{lie group meleshko for Y(FBSDE)} are recovered by solving the Lie's equations
\begin{equation}
    \frac{\partial \Xi}{\partial a}=h(\Xi,\phi,\varphi),\quad  \frac{\partial \varphi}{\partial a}=\gamma(\Xi,\phi,\varphi),\quad \frac{\partial \phi}{\partial a}=\xi(\Xi,\phi,\varphi)
\end{equation}
with initial conditions for $a=0$,
\begin{equation}
    \Xi(t,x,y,0)=t,\quad \varphi(t,x,y,0)=y,\quad \phi(t,x,y,0)=x.
\end{equation}
where $h(t,x, y)=2 \int_{0}^{t} \tau(s,x, y) \mathrm{d} s$ and $
\eta^{2}=\frac{\partial \Xi}{\partial t}.
$

\section{conclusion}
As far as we can tell, applications of the Lie symmetry method to BSDEs remain quite limited and have only recently been explored \cite{ZHANG2021105527,sym11091153}. A promising direction for further research would be to investigate the Lie symmetry of BSDEs with jumps, particularly those of the form 
\begin{equation}
      Y_{t} =\xi+\int_{t}^{T}g(s,Y_s,Z_s,U_s) ds-\int_{t}^{T}Z_s d B_{s}-\int_{t}^{T}\int_{\mathbb{R}}U_s(e)\Tilde{N}(ds,de),\quad t\in [0,T],\label{BSDE with jump}
\end{equation}
where:
\begin{itemize}
    \item 
$(Y_t,Z_t,U_t)$ are the unknown process to be determined.
    \item 
    $B_s$ is a Brownian motion.
    \item 
    $\Tilde{N}(ds,de)$ is the compensated Poisson random measure associated with jumps defined on a probability space.
    \item 
    $g(s,Y_s,Z_s,U_s)$ is the generator of the BSDE, a function that may depend on time $s$, the process $Y_s$, the Brownian motion term $Z_s$, and the jump related-process $U_s(e)$.
    \item 
    $\xi$ is the terminal condition.
\end{itemize}

\noindent If there are no jumps, i.e, if $U_s(e)=0$ then the BSDE reduces to the classical form without jumps:

\begin{equation}
      Y_{t} =\xi+\int_{t}^{T}g(s,Y_s,Z_s) ds-\int_{t}^{T}Z_s d B_{s},\quad t\in [0,T], \label{BSDE without jump}
\end{equation}
Lie symmetries of the BSDE without jump have been studied in \cite{ZHANG2021105527,sym11091153}. However, as far as we know, the symmetries of BSDEs with jumps remains unexplored.\\
More precisely, given a BSDE with jumps \eqref{BSDE with jump}, one may wonder how the compensated Poisson process is affected under the action of a symmetry group (Lie group), or equivalently, under what conditions the transformed compensated Poisson process continues to be a compensated Poisson process under that group action. In the context Brownian motion, Meleshko et al. \cite{srihirun2007definition}, building on \cite{oksendal2013stochastic}, showed that the transformed Brownian motion does indeed stay a Brownian motion (see \ref{new bar(B) Meleshko}).
\section*{Statements and Declarations}

\subsection*{Competing Interests}
The authors declare that they have no competing interests related to this work.

\subsection*{Funding}
This research did not receive any specific grant from funding agencies in the public, commercial, or not-for-profit sectors.

\bibliography{biblio2}

\bibliographystyle{plain}

\end{document}